\documentclass[twoside]{reporteq}
\usepackage{svcon2e}
\usepackage{makeidx}

\usepackage{latexsym}
\usepackage{amsmath}   
\usepackage{amsfonts}
\usepackage{amssymb}
\usepackage{graphicx}
\usepackage{amsthm}

\usepackage{hyperref}
\usepackage[latin1]{inputenc}
\usepackage{url}
\usepackage[numbers]{natbib}

\newcommand{\cl}{C \kern -0.1em \ell}     

\newcommand{\ut}[1]{{\setbox0=\hbox{$#1$}\mathsurround=0pt
       \rlap{\raisebox{-0.8\dp0}{\raisebox{-0.8ex}
       {\kern -0.15ex\hbox{$\tiny\sim$}\kern 0.15ex}}}#1}}
\newcommand{\uti}[1]{{\setbox0=\hbox{$#1$}\mathsurround=0pt
       \rlap{\raisebox{-0.8\dp0}{\raisebox{-0.8ex}
       {\kern -0.3ex\hbox{$\tiny\sim$}\kern 0.3ex}}}#1}}

\newdimen\arrayruleHwidth                 
     \makeatletter
     \def\Hline{\noalign{\ifnum0=`}\fi\hrule \@height \arrayruleHwidth
         \futurelet \@tempa\@xhline}
     \makeatother

\input{psfig}
\newcommand{\ed}{\end{document}}

  \providecommand{\email}[1]{E-mail:
 \href{mailto:#1}{\texttt{#1}}}
 \providecommand{\dprod}{\! \cdot \!}%
 \providecommand{\wprod}{\! \wedge \!}
 
 \providecommand{\pre}[1]{^#1\!}
 \providecommand{\idx}[1]{#1\index{#1}}
\begin{document}

\pagestyle{myheadings} \markboth{José B. Almeida}{Geometric algebra
$G_{4,1}$ and particle dynamics}
\chapter{Geometric algebra and particle dynamics}
\chapterauthors{José B. Almeida}
{ {\renewcommand{\thefootnote}{\fnsymbol{footnote}}
\footnotetext{\kern-15.3pt AMS Subject Classification: 51P05; 81R50;
81Q05.} }}

\begin{abstract}                
In a recent publication \cite{Almeida04:4} it was shown how the
\idx{geometric algebra} $G_{4,1}$\index{G(4,1)@$G_{4,1}$}, the
algebra of \idx{5-dimensional} space-time, can generate relativistic
\idx{dynamics} from the simple principle that only null
\idx{geodesics} should be allowed. The same paper showed also that
\idx{Dirac} equation could be derived from the condition that a
function should be \idx{monogenic} in that algebra; this
construction of the \idx{Dirac} equation allows a choice for the
\idx{imaginary unit} and it was suggested that different
\idx{imaginary units} could be assigned to the various elementary
particles. An earlier paper \cite{Almeida03:5} had already shown the
presence of \idx{standard model} \idx{gauge} group \idx{symmetry} in
complexified ${G}_{1,3}$, an algebra \idx{isomorphic} to
$G_{4,1}$\index{G(4,1)@$G_{4,1}$}.

In this presentation I explore the possible choices for the
\idx{imaginary unit} in the \idx{Dirac} equation to show that
$SU(3)$\index{SU(3)@$SU(3)$} and $SU(2)$\index{SU(2)@$SU(2)$}
symmetries arise naturally from such choices. The \idx{quantum
numbers} derived from the \idx{imaginary unit} are unusual but a
simple conversion allows the derivation of \idx{electric charge} and
\idx{isospin}, \idx{quantum numbers} for two \idx{families} of
particles. This association to elementary particles is not final
because further understanding of the role played by the
\idx{imaginary unit} is needed.

\noindent \textbf{Keywords:} Dirac equation, monogenic functions,
gauge group.
\end{abstract}
%

%
\section{Introduction}
Extracting the laws of physics from a small set of geometric
postulates is an exciting activity which the author has been
exercising for some time with great enthusiasm and led recently to
the presentation of some interesting results \cite{Almeida04:3,
Almeida04:4}. In both of these papers the emphasis is on the
\idx{dynamics} of bodies formulated in \idx{4D} \idx{Euclidean}
space and its relation to general relativity (GR). It was shown how
and when GR and \idx{4D} \idx{Euclidean} metrics can be converted
into each other and an equation was proposed to link the
\idx{Euclidean} metric to its sources. The \idx{4D} \idx{Euclidean}
formulation of \idx{dynamics} is similar to Fermat's principle
applied in 4 dimensions, which motivated the designation
\idx{4-dimensional optics} (\idx{4DO}) to this approach. The point
of departure for the cited work is the 4-dimensional space
constructed with the imposition of \idx{null displacement} in
\idx{5D} space with signature $(-++++).$ Defining \idx{4D} space in
this way leaves its signature undefined until one of the
\idx{coordinates} is expressed in terms of the other 4, allowing for
both \idx{Minkowski} and \idx{Euclidean} signatures and providing
the means for conversion between signatures. The second paper cited
above examined \idx{monogenic} functions \cite{Doran03} in
$G_{4,1}$\index{G(4,1)@$G_{4,1}$} algebra\footnote{I use the
convention that first/second subindex indicates the number of
positive/negative norm \idx{frame} vectors. This convention was
adopted in my previous papers, before I was aware that some people
followed the opposite one.} to show that these verify \idx{Dirac}'s
equation for a free particle. Similarly, if one turns his attention
to \idx{Euclidean} 4-space, the same solutions become plane waves
that extend \idx{4DO} concept into \idx{wave} optics.

It is generally accepted that any physics theory is based on a set
of principles from which predictions are derived using established
mathematical rules; the validity of such theory depends on agreement
between predictions and observed physical reality. In that sense
this paper, as well as those cited above, does not formulate
physical theories because it does not presume any physical
principles. This paper discusses geometry; all along the paper, in
several occasions, a parallel is made with the physical world by
assigning a physical meaning to geometric entities and this allows
predictions to be made. However the validity of derivations and
overall consistency of the exposition is independent of prediction
correctness. The only postulates that are made are of a geometrical
nature and can be summarized in the definition of the space we are
going to work with; this is the \idx{5-dimensional} space with
signature $(-++++),$ where we investigate \idx{monogenic} functions.
Monogenic \index{monogenic} functions provide a more fundamental
constraint than \idx{null subspace} and the latter condition can be
shown to result from the former under certain conditions. The choice
of this geometric space does not imply any assumption for physical
space up to the point where geometric entities like
\idx{coordinates} and \idx{geodesics} start being assigned to
physical quantities like distances and trajectories. Some of those
assignments are carried over from previous work; for instance we
have already established that coordinate $x^0$ is to be taken as
\idx{time} and coordinate $x^4$ as \idx{proper time}
\cite{Almeida04:4}.

Mapping between geometry and physics is facilitated if one chooses
to work always with \idx{non-dimensional} quantities; this is done
with a suitable choice for standards for the fundamental units. From
this point onwards all problems of dimensional homogeneity are
avoided through the use of normalizing factors listed below for all
units, defined with recourse to the fundamental constants: $\hbar
\rightarrow$ Planck constant divided by $2 \pi,$ $G \rightarrow$
\idx{gravitational} constant, $c \rightarrow$ speed of light and $e
\rightarrow$ proton charge.

\begin{center}
\begin{tabular}{c|c|c|c}
Length & Time & Mass & Charge \\
\hline & & & \\

$\displaystyle \sqrt{\frac{G \hbar}{c^3}} $ & $\displaystyle
\sqrt{\frac{G \hbar}{c^5}} $  & $\displaystyle \sqrt{\frac{ \hbar c
}{G}} $  & $e$
\end{tabular}
\end{center}
This normalization defines a system of \emph{\idx{non-dimensional}
units} (\idx{Planck units}) with important consequences, namely: 1)
All the fundamental constants, $\hbar,$ $G,$ $c,$ $e,$ become unity;
2) a particle's Compton frequency, defined by $\nu = mc^2/\hbar,$
becomes equal to the particle's mass; 3) the frequent term
${GM}/({c^2 r})$ is simplified to ${M}/{r}.$

The consideration of \idx{monogenic} functions in the
\idx{5-dimensional} space defined above will be seen to provide
several bridges to the physical world, namely in the areas of
particle \idx{electrodynamics} and the \idx{standard model}. The
exposition makes full use of an extraordinary and little known
mathematical tool called \idx{geometric algebra} (GA), a.k.a.\
\idx{Clifford algebra}, which received an important thrust with the
addition of calculus by David Hestenes \cite{Hestenes84}. A good
introduction to GA can be found in \citet{Gull93} and the following
paragraphs use basically the \idx{notation} and conventions therein
with adaptations to \idx{5D} \idx{hyperbolic space}. A complete
course on physical applications of GA can be downloaded from the
internet \cite{Lasenby99}; the same authors published a more
comprehensive version in book form \cite{Doran03}.
\section{Introduction to \idx{geometric algebra}}
Before we begin our exposition a difficult \idx{notation} issue must
be resolved. We are dealing with \idx{5-dimensional} space but we
are simultaneously interested in two of its 4-dimensional and one
3-dimensional subspaces; ideally our choice of \idx{indices} should
clearly identify their ranges in order to avoid the need to specify
them in every equation. The following diagram shows the index naming
convention used in this paper.

\vspace{12pt} \centerline{\includegraphics[scale=0.8]{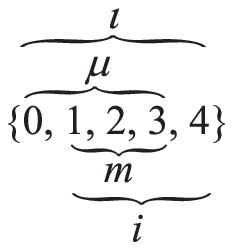}}
\vspace{12pt}

\noindent Indices in the range $\{0,4\}$ will be denoted with Greek
letters $\iota, \kappa, \lambda.$ Indices in the range $\{0,3\}$
will also receive Greek letters but chosen from $\mu, \nu, \xi.$ For
\idx{indices} in the range $\{1,4\}$ we will use Latin letters $i,
j, k$ and finally for \idx{indices} in the range $\{1,3\}$ we will
use also Latin letters chosen from $m, n, o.$

Einstein's summation convention will be adopted as well as the
compact \idx{notation} for partial derivatives $\partial_\iota =
\partial/\partial x^\iota.$ When convenient we will also make the
assignments justified in \citet{Almeida04:4}: $x^0 \equiv t$ and
$x^4 \equiv \tau.$ The \idx{geometric algebra}
$G_{4,1}$\index{G(4,1)@$G_{4,1}$} of the hyperbolic
\idx{5-dimensional} space we want to consider is generated by the
\idx{frame} of \idx{orthonormal} vectors $\sigma_\iota $ verifying
the relations
\begin{align}
    & (\sigma_0)^2  = -1,\\ & \sigma_0 \sigma_i + \sigma_i
    \sigma_0
    =0,\\
    & \sigma_i
    \sigma_j + \sigma_j \sigma_i  = 2 \delta_{i j}.
\end{align}
We will simplify the \idx{notation} for \idx{basis} vector products
using multiple \idx{indices}, i.e.\ $\sigma_\iota \sigma_\kappa
\equiv \sigma_{\iota\kappa}.$ The algebra is 32-dimensional and is
spanned by the \idx{basis}
\begin{itemize}
\item 1 scalar, { $1$}
\item 5 vectors, { $\sigma_\iota$}
\item 10 bivectors (area), { $\sigma_{\iota\kappa}$}
\item 10 trivectors (volume), { $\sigma_{\iota\kappa\lambda}$}
\item 5 tetravectors (4-volume), { $\mathrm{i} \sigma_\iota $}
\item 1 \idx{pseudoscalar} (5-volume), { $\mathrm{i} \equiv \sigma_{01234}$}
\end{itemize}
Several elements of this \idx{basis} square to unity:
\begin{equation}
    (\sigma_i)^2 =  (\sigma_{0i})^2=
    (\sigma_{0i j})^2 =(\mathrm{i}\sigma_0)^2 =1;
\end{equation}
and the remaining square to $-1:$
\begin{equation}
    \label{eq:negative}
    (\sigma_0)^2 = (\sigma_{ij})^2 = (\sigma_{ijk})^2 =
    (\mathrm{i}\sigma_i)^2 = \mathrm{i}^2=-1.
\end{equation}
Note that the \idx{pseudoscalar} $\mathrm{i}$ commutes with all the
other \idx{basis} elements while being a square root of $-1$ and
plays the role of the scalar imaginary in complex algebra. When
writing \idx{matrix} expressions equivalent to \idx{geometric
algebra} we will represent the \idx{complex imaginary} by the symbol
$\mathrm{j},$ following the use of electrical engineering.

The geometric product of any two vectors $a = a^\iota \sigma_\iota$
and $b = b^\kappa \sigma_\kappa$ is evaluated making use of the
distributive property
\begin{equation}
    ab = \left(-a^0 b^0 + \sum_i a^i b^i \right) + \sum_{\iota \neq \kappa}
    a^\iota b^\kappa \sigma_{\iota \kappa};
\end{equation}
and we notice it can be decomposed into a symmetric part, a scalar
called the inner or interior product, and an anti-symmetric part, a
bivector called the outer or exterior product.
\begin{equation}
    ab = a \dprod b + a \wprod b,~~~~ ba = a \dprod b - a \wprod b.
\end{equation}
Reversing the definition one can write inner and outer products as
\begin{equation}
    a \dprod b = \frac{1}{2}\, (ab + ba),~~~~ a \wprod b = \frac{1}{2}\, (ab -
    ba).
\end{equation}
When a vector is operated with a multivector the inner product
reduces the grade of each element by one unit and the outer product
increases the grade by one. By convention the inner product of a
vector and a scalar produces a vector.

The \idx{geometric algebra} $G_{4,1}$\index{G(4,1)@$G_{4,1}$} is
\idx{isomorphic} to the complex algebra of $4 \times 4$
\idx{matrices}, designated by $M(4,{C}).$ When necessary we will
invoke this \idx{isomorphism} to simplify the exposition of some
concepts, namely those related to \idx{symmetry}. The following
mapping between \idx{basis} vectors and \idx{matrices} will be used
\begin{equation}
\label{eq:diracmatrix}
\begin{split}
    \sigma_0 \equiv  \begin{pmatrix} -\mathrm{j} & 0 & 0 & 0 \\
    0 & \mathrm{j} & 0 & 0 \\ 0 & 0 & -\mathrm{j} & 0 \\ 0 & 0 & 0 &
    \mathrm{j}
    \end{pmatrix},~~
    \sigma_1 \equiv & \begin{pmatrix} 0 & 0 & 0 & 1 \\
    0 & 0 & -1 & 0 \\ 0 & -1 & 0 & 0 \\ 1 & 0 & 0 & 0
    \end{pmatrix},~~
    \sigma_2 \equiv   \begin{pmatrix} 0 & \mathrm{j} & 0 & 0 \\
    -\mathrm{j} & 0 & 0 & 0 \\ 0 & 0 & 0 & -\mathrm{j} \\ 0 & 0 & \mathrm{j} & 0
    \end{pmatrix} \\
    \sigma_3 \equiv \begin{pmatrix} 0 & 1 & 0 & 0 \\
    1 & 0 & 0 & 0 \\ 0 & 0 & 0 & 1 \\ 0 & 0 & 1 & 0
    \end{pmatrix} & ,~~
    \sigma_4 \equiv  \begin{pmatrix} 0 & 0 & 0 & -\mathrm{j} \\
    0 & 0 & \mathrm{j} & 0 \\ 0 & -\mathrm{j} & 0 & 0 \\ \mathrm{j} & 0 & 0 & 0
    \end{pmatrix} .
\end{split}
\end{equation}

We will encounter exponentials with multivector exponents; two
particular cases of exponentiation are specially important. If $u$
is such that $u^2 = -1$ and $\theta$ is a scalar
\begin{equation}
\begin{split}
   \mathrm{e}^{u \theta} =& 1 + u \theta -\frac{\theta^2}{2!} - u
    \frac{\theta^3}{3!} + \frac{\theta^4}{4!} + \ldots  \\
    =& 1 - \frac{\theta^2}{2!} +\frac{\theta^4}{4!}- \ldots \{=
    \cos \theta \} \\
    &+ u \theta - u \frac{\theta^3}{3!} + \ldots \{= u \sin
    \theta\}\nonumber \\
    =&  \cos \theta + u \sin \theta.
\end{split}
\end{equation}
Conversely if $h$ is such that $h^2 =1$
\begin{equation}
\begin{split}
    \mathrm{e}^{h \theta} =& 1 + h \theta +\frac{\theta^2}{2!} + h
    \frac{\theta^3}{3!} + \frac{\theta^4}{4!} + \ldots  \\
    =& 1 + \frac{\theta^2}{2!} +\frac{\theta^4}{4!}+ \ldots \{=
    \cosh \theta \}\\
    &+ h \theta + h \frac{\theta^3}{3!} + \ldots \{= h \sinh
    \theta\}  \\
    =&  \cosh \theta + h \sinh \theta.
\end{split}
\end{equation}
The exponential of bivectors is useful for defining rotations; a
rotation of vector $a$ by angle $\theta$ on the $\sigma_{12}$ plane
is performed by
\begin{equation}
    a' = \mathrm{e}^{\sigma_{21} \theta/2} a
    \mathrm{e}^{\sigma_{12} \theta/2}= \tilde{R} a R;
\end{equation}
the tilde denotes reversion and reverses the order of all products.
As a check we make $a = \sigma_1$
\begin{eqnarray}
    \mathrm{e}^{-\sigma_{12} \theta/2} \sigma_1
    \mathrm{e}^{\sigma_{12} \theta/2} &=&
    \left(\cos \frac{\theta}{2} - \sigma_{12}
    \sin \frac{\theta}{2}\right) \sigma_1
    \left(\cos \frac{\theta}{2} + \sigma_{12} \sin
    \frac{\theta}{2}\right)\nonumber \\
    &=& \cos \theta \sigma_1 + \sin \theta \sigma_2.
\end{eqnarray}
Similarly, if we had made $a = \sigma_2,$ the result would have been
$-\sin \theta \sigma_1 + \cos \theta \sigma_2.$

If we use $B$ to represent a bivector belonging to \idx{Euclidean}
space and define its norm by $|B| = (B \tilde{B})^{1/2},$ a general
rotation in 4-space is represented by the \idx{rotor}
\begin{equation}
    R \equiv e^{-B/2} = \cos\left(\frac{|B|}{2}\right) -  \frac{B}{|B|}
    \sin\left(\frac{|B|}{2}\right).
\end{equation}
The rotation angle is $|B|$ and the rotation plane is defined by
$B.$ A \idx{rotor} is defined as a \idx{unitary} even multivector (a
multivector with even grade components only) which squares to unity;
we are particularly interested in rotors with bivector components.
It is more general to define a rotation by a plane (bivector) then
by an axis (vector) because the latter only works in 3D while the
former is applicable in any dimension.

In a general situation the \idx{frame} vectors need not be
orthonormed; we designate such \idx{frame} by \emph{\idx{refractive
index} \idx{frame}} \cite{Almeida04:4} and denote its vectors by the
symbol $g_\iota$
\begin{equation}
    g_\iota = n^\kappa_\iota \sigma_\kappa,
\end{equation}
where $n^\kappa_\iota$ is called the \emph{\idx{refractive index}
tensor}. Refractive index frames are used for \idx{gravitational}
\idx{dynamics}, which lies outside the scope of the present paper.
In association with the \idx{refractive index} \idx{frame} we will
now introduce the concept of \idx{reciprocal frame}, with
\idx{frame} vectors $g^\iota,$ which will be needed in several
occasions along the paper. For our \idx{5-dimensional} space the
\idx{reciprocal frame} is defined by the relations
\begin{equation}
    \label{eq:recframe}
    g^\iota \dprod g_\kappa = {\delta^\iota}_\kappa.
\end{equation}
Following the procedure outlined in \cite{Doran03} to determine the
\idx{reciprocal frame} vectors, we define the 5-volume
\idx{pseudoscalar}
\begin{equation}
    V = \bigwedge_{\iota=0 }^4 g_\iota  = |V|\, \mathrm{i},
\end{equation}
where the big wedge symbol is used to make the exterior product of
the $g_\iota.$ The \idx{reciprocal frame} vectors can then be found
using the formula \cite{Doran03}
\begin{equation}
    g^\iota = \frac{(-1)^{(\iota+1)}}{|V|} \bigwedge_{\iota \neq \kappa} g_\kappa \mathrm{i}.
\end{equation}
Applying the definition to the $\sigma_\iota$ \idx{frame} we are now
working with, it is easy to see that the \idx{reciprocal frame} is
$\sigma^0 = -\sigma_0,$ $\sigma^i = \sigma_i.$

The first use we will make of the \idx{reciprocal frame} is for the
definition of a \idx{vector derivative}. In any space
\begin{equation}
    \mathrm{D} = g^\iota\partial_\iota
\end{equation}
is a vector and as such can be left or right multiplied with other
vectors or multivectors. For instance, when multiplied by vector $a$
on the right the result comprises scalar and bivector terms
$\mathrm{D} a = \mathrm{D} \dprod a + \mathrm{D} \wprod a.$ The
scalar part can be immediately associated with the divergence and
the bivector part is called the exterior derivative; in the
particular case of \idx{Euclidean} 3-dimensional space it is
possible to define the $\mathrm{curl}$ of a vector by
$\mathrm{curl}(a) = -\sigma_{123} \mathrm{D} \wprod a.$

It will be convenient, sometimes, to use vector derivatives in
subspaces of \idx{5D} space; these will be denoted by an upper index
before the $\mathrm{D}$ and the particular index used determines the
subspace to which the derivative applies; For instance $^m\mathrm{D}
= \sigma^m \partial_m = \sigma^1 \partial_1 + \sigma^2 \partial_2 +
\sigma^3 \partial_3.$ In \idx{5-dimensional} space it will be useful
to split the \idx{vector derivative} into its \idx{time} and
4-dimensional parts
\begin{equation}
    \mathrm{D} = -\sigma_0\partial_t + \sigma^i \partial_i
    = -\sigma_0\partial_t
    + \pre{i}\mathrm{D}.
\end{equation}

We define also second order differential operators, generally
designated \idx{Laplacian}, resulting from the product of the
\idx{vector derivative} by itself. The square of a vector is always
a scalar and the \idx{vector derivative} is no exception, so the
\idx{Laplacian} is a scalar operator, which consequently acts
separately in each component of a multivector. For $4+1$ space it is
\begin{equation}
    \mathrm{D}^2 = -\frac{\partial^2}{\partial t^2} + \pre{i}\mathrm{D}^2.
\end{equation}
One sees immediately that a 4-dimensional \idx{wave} equation is
obtained zeroing the \idx{Laplacian} of some multivector function
\begin{equation}
    \label{eq:4dwave}
    \mathrm{D}^2 \psi = \left(-\frac{\partial^2}{\partial t^2} +
    \pre{i}\mathrm{D}^2\right)\psi = 0.
\end{equation}
We will be returning to this \idx{wave} equation but for the moment
we must spend some time examining the symmetries of $G_{4,1}.$
\section{Symmetries of $G_{4,1}$\index{G(4,1)@$G_{4,1}$} algebra}
In this algebra it is possible to find a maximum of four mutually
annihilating \idx{idempotents}, which generate with $0$ an additive
group of order 16; for a demonstration see \citet[section
17.5]{Lounesto01}. Those \idx{idempotents} can be generated by a
choice of two commuting \idx{basis} elements which square to unity;
for the moment we will use $\sigma_{023}$ and $\sigma_{014}.$ The
set of 4 \idx{idempotents} is then given by
\begin{equation}
\begin{split}
    \label{eq:idempotent}
    f_1 = \frac{(1 + \sigma_{023})(1 + \sigma_{014})}{4}, ~~
    f_2 = \frac{(1 + \sigma_{023})(1 - \sigma_{014})}{4}, \\
    f_3 = \frac{(1 - \sigma_{023})(1 - \sigma_{014})}{4}, ~~
    f_4 = \frac{(1 - \sigma_{023})(1 + \sigma_{014})}{4}.
\end{split}
\end{equation}
Using Eq.\ (\ref{eq:diracmatrix}) to make \idx{matrix} replacements
of $\sigma_{023}$ and $\sigma_{014}$ one can find \idx{matrix}
equivalents to these \idx{idempotents}; those are \idx{matrices}
which have only one non-zero element, located on the diagonal and
with unit value.

$SU(3)$\index{SU(3)@$SU(3)$} \idx{symmetry} can now be demonstrated
by construction of the 8 \idx{generators}
\begin{equation}
\label{eq:su3gen}
\begin{split}
    \lambda_1 =& \sigma_{02} (f_1+f_2)= \frac{\sigma_{3}+\sigma_{02}}{2}, \\
    \lambda_2 =& \sigma_{03} (f_1+f_2)= \frac{-\sigma_{2}+\sigma_{03}}{2}, \\
    \lambda_3 =& f_1 - f_2= \frac{\sigma_{014}-\sigma_{1234}}{2}, \\
    \lambda_4 =& -\sigma_1 (f_2 + f_3)= \frac{-\sigma_{1}-\sigma_{04}}{2}, \\
    \lambda_5 =& -\sigma_4 (f_2 + f_3)= \frac{-\sigma_{4}+\sigma_{01}}{2}, \\
    \lambda_6 =& \sigma_{012} (f_1 + f_3)= \frac{\sigma_{012}+\sigma_{034}}{2}, \\
    \lambda_7 =& -\sigma_{024} (f_1 + f_3)= \frac{\sigma_{013}-\sigma_{024}}{2}, \\
    \lambda_8 =& \frac{f_1 + f_2 - 2 f_3}{\sqrt{3}}= \frac{2 \sigma_{023} + \sigma_{014} +
    \sigma_{1234}}{2 \sqrt{3}}.
\end{split}
\end{equation}
These have the following \idx{matrix} equivalents
\begin{equation}
\label{eq:gelmann}
\begin{split}
    \lambda_1 \equiv  \begin{pmatrix} 0 & 1 & 0 & 0 \\
    1 & 0 & 0 & 0 \\ 0 & 0 & 0 & 0 \\ 0 & 0 & 0 & 0
    \end{pmatrix},~~
    \lambda_2 \equiv & \begin{pmatrix} 0 & -\mathrm{j} & 0 & 0 \\
    \mathrm{j} & 0 & 0 & 0 \\ 0 & 0 & 0 & 0 \\ 0 & 0 & 0 & 0
    \end{pmatrix},~~
    \lambda_3 \equiv  \begin{pmatrix} 1 & 0 & 0 & 0 \\
    0 & -1 & 0 & 0 \\ 0 & 0 & 0 & 0 \\ 0 & 0 & 0 & 0
    \end{pmatrix},\\
    \lambda_4 \equiv  \begin{pmatrix} 0 & 0 & 0 & 0 \\
    0 & 0 & 1 & 0 \\ 0 & 1 & 0 & 0 \\ 0 & 0 & 0 & 0
    \end{pmatrix}  ,  ~~
    \lambda_5 \equiv & \begin{pmatrix} 0 & 0 & 0 & 0 \\
    0 & 0 & -\mathrm{j} & 0 \\ 0 & \mathrm{j} & 0 & 0 \\ 0 & 0 & 0 & 0
    \end{pmatrix}, ~~
    \lambda_6 \equiv  \begin{pmatrix} 0 & 0 & 1 & 0 \\
    0 & 0 & 0 & 0 \\ 1 & 0 & 0 & 0 \\ 0 & 0 & 0 & 0
    \end{pmatrix}  \\
    \lambda_7 \equiv  \begin{pmatrix} 0 & 0 & -\mathrm{j} & 0 \\
    0 & 0 & 0 & 0 \\ 0 & \mathrm{j} & 0 & 0 \\ 0 & 0 & 0 & 0
    \end{pmatrix}, & ~~
    \lambda_8 \equiv \left(1/\sqrt{3}\right) \begin{pmatrix} 1 & 0 & 0 & 0 \\
    0 & 1 & 0 & 0 \\ 0 & 0 & -2 & 0 \\ 0 & 0 & 0 & 0
    \end{pmatrix}  ,
\end{split}
\end{equation}
which reproduce \idx{Gell-Mann} \idx{matrices} in the upper-left $3
\times 3$ corner \cite{Almeida03:5, Greiner01, Cottingham98}. Since
the algebra is \idx{isomorphic} to complex $4 \times 4$ \idx{matrix}
algebra, one expects to find higher order symmetries;
\citet{Greiner01} show how one can add 7 additional \idx{generators}
to those of $SU(3)$\index{SU(3)@$SU(3)$} in order to obtain
$SU(4)$\index{SU(4)@$SU(4)$} and the same procedure can be adopted
in \idx{geometric algebra}. We then define the following additional
$SU(4)$\index{SU(4)@$SU(4)$} \idx{generators}
\begin{equation}
\label{eq:su4gen}
\begin{split}
    \lambda_9 =& \sigma_1 (f_1+f_4) = \frac{\sigma_{1}- \sigma_{04}}{2}, \\
    \lambda_{10} =& \sigma_4 (f_1+f_4) = \frac{\sigma_{4}+ \sigma_{01}}{2}, \\
    \lambda_{11} =& -\sigma_{012} (f_2 + f_4) = \frac{-\sigma_{012}- \sigma_{034}}{2}, \\
    \lambda_{12} =& \sigma_{024} (f_2 + f_4) = \frac{\sigma_{013}+ \sigma_{024}}{2}, \\
    \lambda_{13} =& \sigma_{3} (f_3 + f_4) = \frac{\sigma_{3}-\sigma_{02}}{2}, \\
    \lambda_{14} =& \sigma_{2} (f_3 + f_4) = \frac{\sigma_{2}+ \sigma_{03}}{2}, \\
    \lambda_{15} =& \frac{f_1 + f_2 + f_3 - 3 f_4}{\sqrt{6}} =
    \frac{\sigma_{023}
    - \sigma_{014} - \sigma_{1234}}{\sqrt{6}}.
\end{split}
\end{equation}
Once again, making the replacements with Eq.\ (\ref{eq:diracmatrix})
produces the \idx{matrix} equivalent \idx{generators}
\begin{equation}
\label{eq:greiner}
\begin{split}
    \lambda_9 \equiv  \begin{pmatrix} 0 & 0 & 0 & 1 \\
    0 & 0 & 0 & 0 \\ 0 & 0 & 0 & 0 \\ 1 & 0 & 0 & 0
    \end{pmatrix},~~
    \lambda_{10} \equiv & \begin{pmatrix} 0 & 0 & 0 & -\mathrm{j} \\
    0 & 0 & 0 & 0 \\ 0 & 0 & 0 & 0 \\ \mathrm{j} & 0 & 0 & 0
    \end{pmatrix},~~
    \lambda_{11} \equiv  \begin{pmatrix} 0 & 0 & 0 & 0 \\
    0 & 0 & 0 & 1 \\ 0 & 0 & 0 & 0 \\ 0 & 1 & 0 & 0
    \end{pmatrix},\\
    \lambda_{12} \equiv  \begin{pmatrix} 0 & 0 & 0 & 0 \\
    0 & 0 & 0 & -\mathrm{j} \\ 0 & 0 & 0 & 0 \\ 0 & \mathrm{j} & 0 & 0
    \end{pmatrix}  ,  ~~
    \lambda_{13} \equiv & \begin{pmatrix} 0 & 0 & 0 & 0 \\
    0 & 0 & 0 & 0 \\ 0 & 0 & 0 & 1 \\ 0 & 0 & 1 & 0
    \end{pmatrix}, ~~
    \lambda_{14} \equiv  \begin{pmatrix} 0 & 0 & 0 & 0 \\
    0 & 0 & 0 & 0 \\ 0 & 0 & 0 & -\mathrm{j} \\ 0 & 0 & \mathrm{j} & 0
    \end{pmatrix},  \\
    \lambda_{15} \equiv & \left(1/\sqrt{6}\right) \begin{pmatrix} 1 & 0 & 0 & 0 \\
    0 & 1 & 0 & 0 \\ 0 & 0 & 1 & 0 \\ 0 & 0 & 0 & -3
    \end{pmatrix} .
\end{split}
\end{equation}

The \idx{standard model} involves the consideration of two
independent $SU(3)$\index{SU(3)@$SU(3)$} groups, one for
\idx{colour} and the other one for \idx{isospin} and
\idx{strangeness}; if \idx{generators} $\lambda_1$ to $\lambda_8$
apply to one of the $SU(3)$\index{SU(3)@$SU(3)$} groups we can
produce the \idx{generators} of the second group by resorting to the
\idx{basis} elements $\sigma_3$ and $\sigma_{04}.$ The new set of 4
\idx{idempotents} is then given by
\begin{equation}
\begin{split}
    \label{eq:idempotent2}
    f_1 = \frac{(1 + \sigma_3)(1 + \sigma_{04})}{4}, ~~
    f_2 = \frac{(1 + \sigma_3)(1 - \sigma_{04})}{4}, \\
    f_3 = \frac{(1 - \sigma_3)(1 - \sigma_{04})}{4}, ~~
    f_4 = \frac{(1 - \sigma_3)(1 + \sigma_{04})}{4}.
\end{split}
\end{equation}
Again a set of $SU(3)$\index{SU(3)@$SU(3)$} \idx{generators} can be
constructed following a procedure similar to the previous one
\begin{equation}
\label{eq:su3gen2}
\begin{split}
    \alpha_1 =& \sigma_{02} (f_1+f_2)= \frac{\sigma_{02}+\sigma_{023}}{2}, \\
    \alpha_2 =& \sigma_{01} (f_1+f_2)= \frac{\sigma_{01}+\sigma_{013}}{2}, \\
    \alpha_3 =& f_1 - f_2= \frac{\sigma_{04}-\sigma_{034}}{2}, \\
    \alpha_4 =& \sigma_2 (f_2 + f_3)= \frac{\sigma_{2}+\sigma_{024}}{2}, \\
    \alpha_5 =& -\sigma_1 (f_2 + f_3)= \frac{-\sigma_{1}-\sigma_{014}}{2}, \\
    \alpha_6 =& \sigma_4 (f_1 + f_3)= \frac{\sigma_{4}-\sigma_{03}}{2}, \\
    \alpha_7 =& \sigma_{012} (f_1 + f_3)= \frac{\sigma_{012}+\sigma_{1234}}{2}, \\
    \alpha_8 =& \frac{f_1 + f_2 - 2 f_3}{\sqrt{3}}= \frac{2 \sigma_3 + \sigma_{04} +
    \sigma_{034}}{2 \sqrt{3}}.
\end{split}
\end{equation}
This new $SU(3)$\index{SU(3)@$SU(3)$} group is necessarily
independent from the first one because its \idx{matrix}
representation involves \idx{matrices} with all non-zero
rows/columns, while the group generated by $\lambda_1$ to
$\lambda_8$ uses \idx{matrices} with zero fourth row/column. In the
following section we will discuss which of the two groups should be
associated with \idx{colour}.
\section{Monogenic \index{monogenic} functions}
A \idx{monogenic} function \cite{Doran03, Lasenby99} is defined in
any space by the condition $\mathrm{D} \psi = 0$. The
\idx{monogenic} condition is more restrictive then the zero
\idx{Laplacian} condition of Eq.\ (\ref{eq:4dwave}) but obviously a
\idx{monogenic} function always verifies the latter. In particular
we can be sure that a \idx{monogenic} function of $4 + 1$ space
produces a 4-dimensional \idx{wave}, although the reverse may not be
true. Splitting the \idx{vector derivative} into its \idx{time} and
\idx{Euclidean} components, we can write the \idx{monogenic}
condition as
\begin{equation}
    \label{eq:1stwave}
    \mathrm{D}\psi = (\pre{i}\mathrm{D} - \sigma_0\partial_t) \psi = 0.
\end{equation}
We expect plane \idx{wave} solutions and so we try
\begin{equation}
    \label{eq:psidef}
    \psi =   \psi_0 \mathrm{e}^{u (\pm p_0 t + p_i x^i)};
\end{equation}
here $u$ is a square root of $-1$ whose characteristics we shall
determine, $p_0$ is the \idx{wave} \idx{angular frequency} and $p_i$
are components of a generalized \idx{wave vector}. When this
solution is inserted in the 1st order equation (\ref{eq:1stwave}) we
get
\begin{equation}
    \label{eq:verifysol}
    (\sigma^i p_i  \mp \sigma_0 p_0)\psi_0 u = 0.
\end{equation}
The first member can only be zero if $\psi_0$ is a multiple of the
vector in parenthesis and is \idx{nilpotent}, i.e.
\begin{equation}
    \label{eq:nilpotent}
    \sum_i (p_i)^2 -(p_0)^2 =0.
\end{equation}
\idx{Dirac} equation can be recovered from Eq.\ (\ref{eq:1stwave})
if it is multiplied on the left by $\sigma_4$
\begin{equation}
    (-\sigma_{40} \partial_t + \sigma_{4m} \partial_m + \partial_4) \psi=0.
\end{equation}
Now note that $\sigma_{40}$ squares to unity and $\sigma_{4m}$
squares to minus unity, so it is legitimate to make assignments to
the \idx{Dirac matrices}: $\gamma^0 \equiv -\sigma_{40},$ $\gamma^m
\equiv \sigma_{4m}.$ The \idx{matrix} $\gamma^5 = \mathrm{i}
\gamma^0 \gamma^1 \gamma^2 \gamma^3$ becomes, upon substitution,
$\gamma^5 \equiv - \sigma_4.$

The last term in the previous equation must be examined with
consideration for the proposed solution. Deriving $\psi$ with
respect to $x^4$ we get
\begin{equation}
    \label{eq:d4psi}
    \partial_4 \psi
    = p_4 \psi_0 u \mathrm{e}^{u(\pm p_0 t + p_\mu x^\mu)}.
\end{equation}
Since $u$ will always commute with the exponential, this simplifies
to
\begin{equation}
    \label{eq:p4psi}
    \partial_4 \psi = - p_4 \psi u.
\end{equation}
We can then make the further assignments $p_0 = E$ (\idx{energy}),
$p_4 = m$ (\idx{rest mass}), and write
\begin{equation}
    \label{eq:Dirac}
    \gamma^{{\mu}}\partial_{{\mu}} \psi = - m \psi u.
\end{equation}
We note here that  \citet{Rowlands03} has been proposing a
\idx{nilpotent} formulation of \idx{Dirac} equation for some years,
albeit with a different algebra.

The \idx{wave} function in Eq.\ (\ref{eq:psidef}) can now be given a
different form, taking in consideration the previous assignments
\begin{equation}
    \label{eq:fermion}
    \psi = A(\sigma^4 m + \mathbf{p} \mp \sigma_0 E)
     \mathrm{e}^{u (\pm E t + \mathbf{p}\cdot \mathbf{x} + m \tau + \alpha)};
\end{equation}
where $A$ is the amplitude, $\mathbf{p} = \sigma^m p_m $ is the
3-dimensional \idx{momentum} vector, $\mathbf{x} = \sigma_m x^m$ is
the 3-dimensional position and $\alpha$ is the phase angle.

Equation (\ref{eq:1stwave}) was written under the assumption of an
orthonormed \idx{frame} but this need not always be the case; the
consideration of a \idx{refractive index} allows the formulation of
particle \idx{dynamics} under \idx{gravity}.

For \idx{gauge} fields we modify the \idx{vector derivative} by
introducing a \idx{gauge derivative}
\begin{equation}
    \mathcal{D} =\, \pre{\mu}\mathrm{D} + (\sigma^4 + \frac{1}{m} A_\mu \sigma^\mu
    u)\partial_4
\end{equation}

Applying the \idx{gauge derivative} to a general function $\psi$ we
get
\begin{equation}
    \label{eq:dynamics}
    \mathcal{D} \psi = \left(\frac{1}{m} A_\mu
    \sigma^\mu u \partial_4 + \sigma^\iota \partial_\iota \right)
    \psi = 0.
\end{equation}
Notice that the $x^4$ derivative of $\psi$ includes the $u$ factor
which multiplies the same factor present in the \idx{gauge
derivative} operator; when this happens, if additionally $u$
commutes with $\psi$, the parenthesis becomes a vector
\begin{equation}
    g^\iota \partial_\iota \psi = ( -A_\mu
    \sigma^\mu  + \sigma^\iota \partial_\iota )
    \psi = 0.
\end{equation}
Multiplying on the left and replacing with \idx{gamma matrices}, as
above, we get \idx{Dirac} equation in an electromagnetic field
(remember the electron has charge $-1)$.
\begin{equation}
    \gamma^\mu (\partial_\mu -  A_\mu) \psi = - m \psi u.
\end{equation}
Obviously the parenthesis in Eq. (\ref{eq:dynamics}) only becomes a
vector if the imaginary factor $u$ in the \idx{gauge derivative}
operator is the same as in the wavefunction. When this does not
happen, the remaining non-vector term does not play a role in
\idx{dynamics}.

We must now examine in detail what sort of geometrical elements can
be used as imaginary $u$ in the solution (\ref{eq:psidef}); it will
be more convenient to write the imaginary as a product of the
\idx{pseudoscalar} and a \idx{unitary} element, $u = \mathrm{i} h,$
since we have already established that each set of 4
\idx{idempotents} is generated by a pair of commuting \idx{unitary}
\idx{basis} elements. Let any two such \idx{basis} elements be
denoted as $h_1$ and $h_2;$ then the product $h_3 = h_1 h_2$ is
itself a third commuting \idx{basis} element. For consistence we
choose, as before,
\begin{equation}
    h_1 \equiv \sigma_{023},~~~~h_2\equiv \sigma_{014};
\end{equation}
to get
\begin{equation}
    h_3 \equiv \sigma_{1234},
\end{equation}
which commutes with the other two as can be easily verified. The
result of this exercise is the existence of \idx{triads} of
commuting \idx{unitary} \idx{basis} elements but no tetrads of such
elements. We are led to state that a general \idx{unitary} element
is a linear combination of unity and the three elements of one
\idx{triad}
\begin{equation}
\label{eq:unitary}
    h = a_0 + a_1 h_1 + a_2 h_2 + a_3 h_3.
\end{equation}
Since $h$ is \idx{unitary} and the three $h_m$ commute we can write
\begin{equation}
\label{eq:uniteq}
\begin{split}
    h^2 =\,& \left[(a_0)^2 + (a_1)^2 + (a_2)^2 + (a_3)^2\right] + 2 (a_0 a_1 - a_2
    a_3) h_1\\
    & + 2 (a_0 a_2 - a_1 a_3) h_2 + 2 (a_0 a_3 - a_1 a_2) h_3
    =1
\end{split}
\end{equation}
The only form this equation can be verified is if the term in square
brackets is unity while all the others are zero. We then get a set
of four simultaneous equations with a total of sixteen solutions, as
follows: 8 solutions with one of the $a_\mu$ equal to $\pm 1$ and
all the others zero, 6 solutions with two of the $a_\mu$ equal to
$-1/2$ and the other two equal to $1/2$ and 2 solutions with all the
$a_\mu$ simultaneously $\pm 1/2.$ The $a_\mu$ coefficients play the
role of \idx{quantum numbers} which determine the particular
\idx{imaginary unit} that goes into Eq.\ (\ref{eq:psidef}) and
consequently into \idx{Dirac}'s equation; these unusual \idx{quantum
numbers} can be combined in order to produce a more conventional set
of numbers, as we shall see next. In \idx{matrix} form these
solutions correspond to \idx{unitary} diagonal \idx{matrices} with
0, 1, 2, 3 or 4 negative elements.

There are 16 possible \idx{unitary} elements defined by Eq.\
(\ref{eq:unitary}), corresponding to different solutions of Eq.\
(\ref{eq:uniteq}). If desired it is possible to express these
solutions in terms of the $SU(4)$\index{SU(4)@$SU(4)$}
\idx{generators} $\lambda_3,$ $\lambda_8$ and $\lambda_{15}$ in
order to highlight the symmetries. The conversion relations for the
generator coefficients are: $\lambda_3 \rightarrow a_2 - a_3,$
$\lambda_8 \rightarrow (2 a_1 + a_2 + a_3)/\sqrt{3}$ and
$\lambda_{15} \rightarrow \sqrt{2/3}(a_1 -a_2 - a_3).$ In table
\ref{t:quantum} we list all the 16 \idx{unitary} elements and their
respective coefficients.

\begin{table}[htb]
{\tiny
 \caption{\label{t:quantum} Coefficients for the various
\idx{unitary} elements.}
\begin{center}
\begin{tabular}{ r |r r r|r r r| r r| c }
  \hline
   $1$ & $\sigma_{023}$  & $\sigma_{014}$   &   $\sigma_{1234}$  & $\lambda_3$
     & $\lambda_8$  & $\lambda_{15}$  & $q$ & $i_3$ &  \\
   $(a_0)$ & $(a_1)$  &  $(a_2)$  &  $(a_3)$ &  & &
    &
    &  & Designation \\
  \hline
   $1$ & $0$  & $0$  & $0$  & $0$  & $0$  & $0$  & $0$ & $0.5$ &  \\
   $0$ & $1$  & $0$  & $0$  & $0$  & $2/\sqrt{3}$  & $\sqrt{2/3}$  & $2/3$ & $0.5$ & up  \\
   $0$ & $0$  & $1$  & $0$  & $1$  & $1/\sqrt{3}$  & $-\sqrt{2/3}$  & $1/3$ & $0.5$ & anti-down  \\
   $0$ & $0$  & $0$  & $1$  & $-1$  & $1/\sqrt{3}$  & $-\sqrt{2/3}$  &  $1$ & $0.5$ & positron  \\
   $-1$ & $0$  & $0$  & $0$  & $0$  & $0$  & $0$  & $0$ &  $-0.5$ & \\
   $0$ & $-1$  & $0$  & $0$  & $0$  & $-2/\sqrt{3}$  & $-\sqrt{2/3}$  & $-2/3$ & $-0.5$ & anti-up  \\
   $0$ & $0$  & $-1$  & $0$  & $-1$  & $-1/\sqrt{3}$  & $\sqrt{2/3}$  & $-1/3$ & $-0.5$ & down  \\
   $0$ & $0$  & $0$  & $-1$  & $1$  & $-1/\sqrt{3}$  & $\sqrt{2/3}$  & $-1$ & $-0.5$ & electron  \\
   $-1/2$ & $-1/2$  & $1/2$  & $1/2$  & $0$  & $0$  & $-\sqrt{3/2}$  & $1/3$ & $0$ & anti-strange \\
   $-1/2$ & $1/2$  & $-1/2$  & $1/2$  & $-1$  & $1/\sqrt{3}$  & $1/\sqrt{6}$  & $2/3$ & $0$ & charm  \\
   $-1/2$ & $1/2$  & $1/2$  & $-1/2$  & $1$  & $1/\sqrt{3}$  & $1/\sqrt{6}$  & $0$ & $0$ &   \\
   $1/2$ & $-1/2$  & $-1/2$  & $1/2$  & $-1$  & $-1/\sqrt{3}$  & $-1/\sqrt{6}$  & $0$ &  $0$ &  \\
   $1/2$ & $-1/2$  & $1/2$  & $-1/2$  & $1$  & $-1/\sqrt{3}$  & $-1/\sqrt{6}$  & $-2/3$ & $0$ & anti-charm  \\
   $1/2$ & $1/2$  & $-1/2$  & $-1/2$  & $0$  & $0$  & $\sqrt{3/2}$  & $-1/3$ & $0$ & strange \\
   $1/2$ & $1/2$  & $1/2$  & $1/2$  & $0$  & $2/\sqrt{3}$  & $-1/\sqrt{6}$  & $1$ & $1$ & anti-mu  \\
   $-1/2$ & $-1/2$  & $-1/2$  & $-1/2$  & $0$  & $-2/\sqrt{3}$  & $1/\sqrt{6}$  & $-1$ & $-1$ & mu  \\
  \hline
\end{tabular}
\end{center}
    }
\end{table}

In columns labelled $q$ and $i_3$ we make a tentative association
between the coefficients $a_\mu$, \idx{electric charge} and
\idx{isospin} through the formulas $q = 2 a_1 + a_2 + 3 a_3$ and
$i_3 = (a_0 + a_1 + a_2 + a_3)/2$. This association is rather
speculative and likely to change in the future. Electric charge was
assumed to be a count of the number of spatial \idx{basis} vectors
present in the \idx{unitary} element $h$ and \idx{isospin} was
defined in such a way that it would allow the right combinations for
all the particles in two \idx{families}. Further insight onto the
role of the \idx{unitary} element is needed in order to properly
justify the correct association to \idx{quantum numbers}. The
provisional particles' names are included in the last column and we
left four positions blank because their significance is unclear for
the moment.

\section{Conclusion}
Using the algebra of \idx{5-dimensional} spacetime $G_{4,1}$ as a
point of departure, it was possible to derive a number of results of
great physical significance from the properties of that algebra. The
\idx{monogenic} condition applied to functions was seen to allow
\idx{Euclidean} space 4-dimensional plane \idx{wave} solutions but
the same condition could be converted into \idx{Dirac} equation,
under certain circumstances, providing an \idx{Euclidean}
interpretation of \idx{Dirac} spinors. Particle \idx{dynamics} could
be modelled through the definition of a \idx{gauge derivative}
operator, through which the electromagnetic case was successfully
formulated.

The exponent in the solutions of the \idx{monogenic} condition
includes a square root of $-1$ which can be chosen among different
possibilities provided by $G_{4,1}$ algebra. An exploration of those
possibilities showed that they exhibit the same symmetries as the
\idx{standard model} \idx{gauge} group, which allowed a tentative
mapping to the elementary particles of the two lower \idx{families}.


%

  \bibliographystyle{unsrtbda}   
  \bibliography{Abrev,aberrations,assistentes}   
 \small
 \vskip 1pc
 {\obeylines
 \noindent José B. Almeida
 \noindent Universidade do Minho
 \noindent Physics Department
 \noindent 4710-057 Braga
 \noindent Portugal
 \noindent \email{bda@fisica.uminho.pt}
 \vskip 1pc }
 \vskip 6pt
 \noindent Submitted: TBA; Revised: \today.\\

 \normalsize \noindent
\printindex
\end{document}